\documentclass[12pt,leqno]{article}
\usepackage{amsmath,amsthm}
\usepackage{amsfonts,amssymb}
\newtheorem{theorem}{Theorem}
\newtheorem{lemma}[theorem]{Lemma}
\newtheorem{proposition}[theorem]{Proposition}
\newtheorem{corollary}[theorem]{Corollary}

\numberwithin{equation}{section}

\begin{document}
\begin{center}
{\LARGE Dynamics and  Endogeny for recursive processes \\ on trees}\\
\vspace{.1in}
{\large JON WARREN}\\
\vspace{.1in} Department of Statistics, University of Warwick, Coventry CV4 7AL, UK.
\end{center}

\begin{abstract}
We consider a stochastic process $\bigl(\xi_u; u \in \Gamma_\infty \bigr)$ where $\Gamma_\infty$ is the set  of vertices of an infinite 
binary tree 
 which satisfies the recursion relation  
\[
\xi_u= \phi(\xi_{u0},\xi_{u1}, \epsilon_u) \text { for each } u \in \Gamma_\infty.
\]
Here $u0$ and $u1$ denote the two immediate daughters of the vertex  $u$. The random variables $\bigl( \epsilon_u; u\in 
\Gamma_\infty\bigr)$, which 
are
 to be thought of as innovations, are supposed  independent and identically distributed.  A recent paper of Aldous and  Bandyopadhyay has 
drawn 
attention  to the issue of
 endogeny: that is whether the process $\bigl( \xi_u; u \in \Gamma_\infty\bigr)$ is measurable with respect to the innovations process. 
Here we
 restrict attention to the case where each $\xi_u$ takes values in a finite set $S$ and show how
 this question is related to the existence of certain  dynamics.  Using this we develop a necessary and sufficient condition for endogeny 
in terms of 
the coupling rate for a Markov chain on $S^2$ for which the diagonal is absorbing. 
 \end{abstract}
 
\section{Introduction}

We suppose that we are given  probability spaces
$\bigl(S, {\cal S}, \mu \bigr)$ and  $\bigl(E, {\cal E}, \nu\bigr)$ and a mapping
\begin{equation}
\phi: S \times S \times E \rightarrow S
\end{equation}
under which the image of  the product measure $\mu \otimes \mu \otimes \nu $  is $\mu$. 
In general  $\phi: S^n \times E 
\rightarrow S$ may be considered, even with $n$ being infinity, but here we take $n=2$ for simplicity and we further assume that $\phi$ is 
symmetric in its 
first two 
arguments. 
Typically  in applications 
$\nu$ and $\phi$ are given and the measure $\mu$ is to be determined. $\mu$ is then said to be a solution to a  recursive distributional 
equation. Existence 
and uniqueness of $\mu$ is well studied in a number of important cases; the survey paper by Aldous and Bandyopdhyay, \cite{ab}, describes 
many examples.  Here our 
concerns are somewhat different. Associated with this setup is a recursive tree process constructed as follows. Let $\Gamma_\infty= 
\cup_{k\geq 0} 
\{0,1\}^k$ be the set of vertices  of  the  infinite rooted binary tree; 
in particular  $\emptyset\in \Gamma_\infty$ denotes the root. The set of vertices on level $n$ of the tree, $\{0,1\}^n$, is denoted  by 
$G_n$, and 
$\Gamma_n=\cup_{k=0}^n \{0,1\}^k$ denotes the set of
vertices belonging to levels up to and including $n$. A vertex $u=(u_1,u_2,\ldots, u_n)$  has two daughters: vertices 
$u0=(u_1,\ldots,u_n,0)$ and
$u1=(u_1,\ldots,u_n,1)$.  Let $\bigl(\xi_u; u\in \Gamma_\infty \bigr)$ and $\bigl(\epsilon_u; u \in \Gamma_\infty\bigr)$ be the 
co-ordinate maps on  
$S^{\Gamma_\infty} \times E^{\Gamma_\infty}$. Define
\begin{equation}
\Omega= \bigl\{ \omega \in S^{\Gamma_\infty} \times E^{\Gamma_\infty}: \xi_u(\omega)=\phi(\xi_{u0}(\omega),\xi_{u1}(\omega), 
\epsilon_u(\omega)) \text { 
for each } u \in \Gamma_\infty \bigr\},
\end{equation}
and let ${\cal F}$ denote the  restriction of the product $\sigma$-algebra to $\Omega$. There is a unique  probability measure $m$ on 
$\bigl(\Omega,{\cal F}\bigr)$ under which for each $n$,  the joint law of $\bigl(\xi_u; u \in G_n\bigr)$ and $\bigl(\epsilon_u; u \in 
\Gamma_{n-1}\bigr)$ is $\otimes^{G_n} \mu \otimes^{\Gamma_{n-1}} \nu$.  We think of independent random variables $\bigl(\epsilon_u; u \in 
\Gamma_\infty\bigr)$, sometimes referred to as the innovations,  as the input that drives the system, while $\bigl(\xi_u; u \in 
\Gamma_\infty\bigr)$ is thought of  as the response. Aldous and Bandyopadhyay draw attention to investigating  when the  recursive 
tree process is 
endogenous, that is to say every random variable on $\bigl(\Omega,{\cal F}\bigr)$ is $m$-almost surely equal to a function of the 
innovations process 
$\bigl(\epsilon_u; u \in \Gamma_\infty \bigr)$ alone.
Loosely speaking this means there is no additional randomness in the system ``located'' at the boundary of the infinite tree. However care 
must be taken in interpreting this: the tail $\sigma$-algebra
\begin{equation}
\bigcap_{n} \sigma\bigl( \xi_u,\epsilon_u; u \in \Gamma_\infty \setminus \Gamma_n \bigr) 
\end{equation}
is typically  empty even if endogeny does not hold! 
A strong parallel may be drawn  with certain stochastic differential equations which 
admit weak but not strong solutions, see in particular Tsirelson's example discussed in Section V.18 of \cite{rw}.
In  non-endogenous cases it is natural to try to give some explicit description of the additional randomness,  however this does not 
seem possible in any generality. One  case that is  amenable to such an analysis is the linear ``smoothing transformation'' which has been 
extensively studied. Various representations theorems obtained by Durrett and Liggett, \cite{dur}, Liu \cite{liu},  Caliebe and R\"{o}sler 
\cite{cal}, and others reveal a structure in which the innovations are augmented with   Gaussian or Poisson noise living on the boundary of 
the tree.  Another case whose structure can similarly be understood is studied by Biggins \cite{big}.

Let $u_0=\emptyset, u_1,u_2,\ldots, u_n,\ldots$ be an infinite sequence of vertices,  $u_{n+1}$ being a daughter of $u_n$ for each $n$. 
For $n\leq 0$ 
let $\xi_n=\xi_{u_{-n}}$. The law of the sequence $\bigl( \xi_{n}; n \leq 0\bigr)$, which by the symmetry of $\phi$ does not depend on 
the choice of 
sequence of vertices,  is easily seen to be that of a stationary Markov chain (indexed by negative time) with transition kernel $P$ 
defined on  $S\times {\cal 
S}$ 
given by
\begin{equation}
P(x_0,A)= \int_{S} \int_E {\mathbf 1}\bigl( \phi(x_0, x_1,z)\in A\bigr)  \nu(dz)\mu(dx_1).
\end{equation}
In this paper we give a criterion for endogeny is in terms of the corresponding ``two point motion'', that is to say a Markov chain on $S^2$ 
with transition 
kernel 
$P^{(2)}$ on $S^2 \times {\cal S}^2$ given by
\begin{equation}
\label{two}
P^{(2)}((x_0,x_0^\prime),A\times A^\prime)= \int_{S} \int_E {\mathbf 1}\bigl( \phi(x_0, x_1,z)\in A , \phi(x_0^\prime, x_1,z)\in A^\prime 
\bigr) 
\nu(dz)\mu(dx_1).
\end{equation} 
Let $S^\nearrow$ be the diagonal of $S^2$ which is an absorbing set for $P^{(2)}$. Let $P^{(-)}$ be the restriction of $P^{(2)}$ to 
$S^2\setminus 
S^\nearrow$.  
\begin{center}
{\em We now assume that $S$ is finite. }
\end{center}
In this case we can identify $P^{(-)}$ with a non-negative  square matrix to which we can  apply Perron-Frobenius theory. Let $\rho$ be 
the largest 
eigenvalue of $P^{(-)}$.
\begin{theorem} 
\label{main} If 
$2\rho < 1$  then the tree process is endogenous, whereas
if $2 \rho>1$ then the  process is not endogenous. In the critical case $2\rho=1$, if the  additional  conditions \eqref{nondegen1} and 
\eqref{nondegen2} hold, then the tree process is endogenous.
\end{theorem}

Aldous and  Bandyopadhyay have given a  different necessary and sufficient condition for endogeny, the above result being related to 
their condition by 
linearization.  We now elucidate this connection. Associated with $\phi:S\times S \times E \rightarrow S$ is the  ``two-point map'' 
$\phi^{(2)}:S^2\times S^2 \times E \rightarrow S^2$ given by 
\begin{equation}
\phi^{(2)}((x_0,x_0^\prime),(x_1,x_1^\prime),z)= \bigl( \phi(x_0,x_1,z),\phi(x_0^\prime,x_1^\prime,z)\bigr).  
\end{equation}
Now we can define $T^{(2)}:{\cal M}( S^2) \times {\cal M}(S^2) \rightarrow {\cal M}(S^2)$ by setting  $T^{(2)}(\lambda_0,\lambda_1)$  
equal to the image of 
the measure 
$\lambda_0\otimes\lambda_1\otimes \nu$ under the map
$\phi^{(2)}$. Recall that $\mu$ denotes our given measure on $S$ that is invariant for the recursive distributional equation 
corresponding to $\phi$ and 
$\nu$. The  measure  $\mu^\nearrow$ on $S^2$, carried by the diagonal and having marginals $\mu$, is a fixed 
point of $T^{(2)}$ in that $T^{(2)}(\mu^\nearrow, \mu^\nearrow)= \mu^\nearrow$. 
According to Theorem 11 of \cite{ab}, under some minor technical condition, $\mu^\nearrow$ being the only probability measure on 
$S^2$ that is a 
fixed point of $T^{(2)}$, and whose marginals are both equal to $\mu$, is a necessary and sufficient condition for endogeny.
Bandyopadhyay, \cite{band}, has used this criterion to show the endogeny of  certain important examples that arise from applications.
The map $T^{(2)}$  is bilinear and symmetric so that we have for real $\delta$ and $\lambda \in {\cal M}(S^2)$ 
\begin{equation}
T^{(2)}(\mu^\nearrow+\delta \lambda,\mu^\nearrow+\delta \lambda)= \mu^\nearrow +2\delta T^{(2)}( \lambda,\mu^\nearrow)+ \delta^2 
T^{(2)}(\lambda,\lambda),
\end{equation}
which shows that $2T^{(2)}( \cdot,\mu^\nearrow): {\cal M}(S^2)\rightarrow {\cal M}(S^2)$ is the derivative of 
$\lambda \mapsto T^{(2)}(\lambda,\lambda)$ at the fixed 
point $\mu^\nearrow$. 
It is straightforward to verify that
\begin{equation}
\lambda P^{(2)}= T^{(2)}( \lambda,\mu^\nearrow).
\end{equation}
Thus the criterion given in Theorem \ref{main} can be viewed in terms of the stability of the fixed point $\mu^\nearrow$ of $T^{(2)}$, 
and it is quite 
natural that the existence of other fixed points of $T^{(2)}$ should be related to this.

 Let ${\cal H}= L^2(\Omega,{\cal F},m)$, and ${\cal K} \subseteq {\cal H}$ be the subspace  of (equivalence classes of) random variables 
measurable with 
respect to 
the innovations process $\bigl( \epsilon_u; u \in \Gamma_\infty\bigr)$. Endogeny means ${\cal H}={\cal K}$. We also have subspaces ${\cal 
K}_n$ and ${\cal 
H}_n$ containing random variables measurable with respect to $\bigl( \epsilon_u; u \in \Gamma_{n-1}\bigr)$ and with respect to $\bigl( 
\epsilon_u; u \in 
\Gamma_{n-1}\bigr)$ together with $\bigl( \xi_u; u\in G_n\bigr)$ respectively.

In the following three sections we treat the three cases: subcritical $2\rho<1$; supercritical $2 \rho>1$; and critical $2 \rho=1$ 
separately. The subcritical 
case
is treated by using an operator version of Markov's inequality to bound from above the quantity $\|f\|^2-\|P_{{\cal K}_n}f\|^2$ where 
$P_{{\cal K}_n}f$ is 
the orthogonal projection of $f$ onto ${\cal K}_n$.   
The supercritical case is treated by using the left eigenvector of $P^{(-)}$ corresponding to $\rho$ to construct a consistent family of 
quadratic forms on ${\cal H}$, and hence   an operator $Q_\infty$. This operator  is the generator of dynamics that  fix 
the innovations while perturbing the additional randomness at the boundary. The existence of such dynamics precludes endogeny.

In the final section of the paper we show that the dynamics associated with $Q_\infty$ naturally arise through a passage to the limit. 
Consider for each $n \geq 1$  the finite configuration space
 \begin{equation}
\Omega_n= \bigl\{ \omega \in S^{\Gamma_n} \times E^{\Gamma_{n-1}}: \xi_u(\omega)=\phi(\xi_{u0}(\omega),\xi_{u1}(\omega), 
\epsilon_u(\omega)) \text { for each } u \in \Gamma_{n-1} \bigr\},
\end{equation}
 which we equip with the measure $m_n$ under which the joint law of $\bigl( \xi_u; u \in G_n \bigr)$ and $\bigl(\epsilon_u; u \in 
\Gamma_{n-1}\bigr)$ is $\otimes^{G_n} \mu \otimes^{\Gamma_{n-1}} \nu$. Consider the dynamics on $\Omega_n$ with invariant measure $m_n$ 
under which the co-ordinates $\xi_u$ for $u \in G_n$ are independently refreshed at rate $1$, whilst $\xi_v$ for $v \in \Gamma_{n-1}$ are 
determined by application of the map $\phi$ using the innovations $(\epsilon_u; u \in \Gamma_{n-1})$ which are held fixed for all time.
Let $-A_n$ denote the corresponding  generator acting on $L^2(\Omega_n,m_n)$. We show that at the level of generators, 
these dynamics, when slowed down by a factor of $(2\rho)^n$, converge to those associated with $Q_\infty$. 
Let $P_{{\cal H}_n}:{\cal H} \rightarrow {\cal H}_n$ be the orthogonal projection onto ${\cal H}_n$ which we identify with  
$L^2(\Omega_n,m_n)$.  
\begin{theorem}
\label{thmcon}
Suppose that $P^{(-)}$ is primitive, and that $2\rho>1$, then
\[
(2\rho)^{-n} A_n P_{{\cal H}_n} \rightarrow -Q_\infty,
\]
in the strong resolvent sense as $n$ tends to infinity.
\end{theorem}

There is strong similarity between the methods used in this paper
 and those used to study  multitype branching processes. To make this connection fix $f \in {\cal H}_0$  that 
satisfies $\|f\|^2=1$ and 
consider the probability distribution distributions $\mu^{(n)}_f$ on the non-negative integers  defined by
\begin{equation}
\label{specmeasure}
\mu_f^{(n)}(k)=\bigl\| P_{E_k} f \bigr\|^2 \text{ for } k \geq 0,
\end{equation}
where $P_{E_k}$ is the projection operator associated with the  eigenspace $E_k$ of the operator $A_n$ corresponding to the eigenvalue $k$. 
See below at \eqref{number} for further information on this. We use  arguments that  are based on treating this distribution as if it   were 
that of the number of particles alive in generation $n$ of a branching process. In particular Theorem \ref{main} corresponds exactly to the 
fact
 the certainty or otherwise of eventual extinction for a  branching process depends upon 
whether the Malthusian parameter of the process is greater than one.  
The arguments used are also  closely related  to the spectral methods employed by Tsirelson, \cite{t}, to study continuous products of 
probability spaces. In fact  Tsirelson and Vershik  introduced certain  recursive tree processes in \cite{tv} as examples in the theory 
of continuous 
products.

\section{The subcritical case: endogeny}

Recall the probability measure $m$  on $\bigl(\Omega,{\cal F}\bigr)$ is characterized by the fact that the joint law of $\bigl( \xi_u;  
u\in G_n\bigr)$ and  
$\bigl( \epsilon_u; u \in \Gamma_{n-1}\bigr)$ is $\otimes^{G_n} \mu \otimes^{\Gamma_{n-1}} \nu$.  Thus  there are isomorphisms 
between 
Hilbert spaces:
\begin{multline}
\label{iso}
{\cal H}_n \cong L^2\bigl( S^{G_n}\times E^{\Gamma_{n-1}},{\cal S}^{G_n} \times {\cal E}^{\Gamma_{n-1}}, \otimes^{G_n} 
\mu\otimes^{\Gamma_{n-1}} \nu\bigr) \\
 \cong \sideset{}{^{G_n}}\bigotimes L^2(S,{\cal S},\mu) \sideset{}{^{\Gamma_{n-1}}}\bigotimes L^2(E,{\cal E},\nu).
\end{multline}
$L^2(S,{\cal S},\mu)$ is the direct sum of the one-dimensional subspace  of constants together with its orthogonal complement to be 
denoted by $L^2_0(S,{\cal 
S},\mu)$. Decomposing each copy of  $L^2(S,{\cal S},\mu)$ appearing on the righthandside of \eqref{iso} in this manner we obtain
\begin{equation}
\label{subsetdec}
{\cal H}_n= \bigoplus_{S \subseteq G_n} {\cal H}_{S},
\end{equation}
where  for each subset  $S\subseteq G_n$  of vertices on level $n$ of the tree  the corresponding  subspace   ${\cal H}_S$ of the Hilbert 
space ${\cal H}_n$  
is generated by vectors of the form $ \otimes_{u\in G_ n} f_u \otimes_{v\in \Gamma_{n-1} }g_v$ with  $f_u\in L^2_0(S,{\cal S},\mu)$ for 
$u \in S$,  $f_u$ 
being a constant vector  for $u \not\in S$, and $g_v \in L^2(E,{\cal E},\nu)$ being unrestricted.

Given a  linear operator $L$ acting on $L^2(S, {\cal S}, \mu)$  and a vertex $u\in G_n$  we may
consider  ``$L$ applied at $u$''. More precisely we define
an operator $L^{(u)}$ acting  on ${\cal H}_n$  as being unitary equivalent, via the isomorphism \eqref{iso}, to  the tensor product of 
$L$ acting on the copy 
of  $L^2(S, {\cal S}, \mu)$ 
corresponding to $u$ together with the identity on all other factors.  We  consider the case that $L$ is given by $P_{{\mathbf 1}^\bot}$ 
which is the 
orthogonal projection  onto the subspace $L^2_0(S,{\cal S},\mu)$. 
We define an operator on ${\cal H}_n$ via
\begin{equation}
\label{number}
A_n= \sum_{u\in G_n}  P_{{\mathbf 1}^\bot}^{(u)}.
\end{equation}
By considering its action on the generating vectors for each subspace ${\cal H}_S$ we find that the eigenvalues of $A_n$ are 
$0,1,2,\ldots ,2^n$, with the 
eigenvalue $k$ having corresponding eigenspace $\bigoplus_{|S|=k} {\cal H}_S$. Here $|S|$ denotes the number of vertices belonging to 
$S$. 
On the other hand the operator $I-P_{{\cal K}_n}$, which commutes with $A_n$, has eigenvalues $0$ and $1$ with corresponding 
eigenspaces ${\cal K}_n={\cal 
H}_\emptyset$ and $\bigoplus_{S \neq \emptyset} {\cal H}_S$. By decomposing $f \in {\cal H}_n$ according to the subspaces ${\cal H}_S$  
we deduce the 
following inequality. 
\begin{proposition}
\label{inequality}
For every $f \in {\cal H}_n$,
\[
0\leq \|(I- P_{{\cal K}_n})f\|^2 \leq \bigl( f,A_nf\bigr).
\]
\end{proposition}
 
\begin{proof}[Proof of Theorem 1: Subcritical case.]
Because of the recursive structure endogeny is is equivalent to ${\cal H}_0\subset {\cal K}$. 
So fix $f \in {\cal H}_0$, which, in a slight abuse of notation, we also treat as  an element of $L^2(S,{\cal S},\mu)$. To prove $f$ belongs 
to ${\cal K}$ it is 
enough, by virtue of 
Proposition \ref{inequality}, to prove that
\[
\bigl( f ,A_n f  \bigr) \rightarrow  0 \text { as } n \rightarrow \infty.
\]
 We may  express $\bigl( f  , A_n f   \bigr)$ using
couplings  of the tree-indexed process. Fix a vertex $u_n \in G_n$.
Let $\Omega_n^\prime$ be a copy of the finite configuration space $\Omega_n$ and consider the product space $\Omega_n \times 
\Omega^\prime_n$ equipped with co-ordinate maps $(\epsilon_u, \epsilon^\prime_u; u \in \Gamma_{n-1})$ and $(\xi_u,\xi_u^\prime; u \in 
\Gamma_n)$. Let the probability measure $\tilde{m}_n$ on this  product space have both marginals equal to $m_n$ and be such that
\begin{description}
\item $\tilde{m}_n$ is supported on the set where  $\epsilon_u=\epsilon^\prime_u$ for all $u \in \Gamma_{n-1}$ and $\xi_u=\xi^\prime_u$ for 
all $u \in G_n$ except $u_n$;
\item  $( \epsilon_u; u\in \Gamma_{n-1})$, $(\xi_u; u \in \Gamma_n)$ and  $\xi^\prime_{u_n}$ are independent under $\tilde{m}_n$. 
\end{description}
It is easily verified that
\[
\bigl( f , P_{{\mathbf 1}^\bot}^{(u_n)}f )=\tfrac{1}{2} \int_{\Omega_n \times \Omega_n^\prime} (f \circ\xi_\emptyset-f\circ 
\xi_\emptyset^\prime)^2 \;d \tilde{m}_n=(\mu \otimes \mu)P_n^{(2)}g,
\]
where $P_n^{(2)}$ is the $n$-step transition matrix for $P^{(2)}$ and $g$ is the function $g(x,x^\prime)= \tfrac{1}{2}(f(x)-f(x^\prime))^2$ 
on $S^2$. 
Summing over the possible choices of $u_n$ we obtain
\[ 
\bigl( f ,A_n f  \bigr)=2^{n}(\mu \otimes \mu)P_n^{(2)}g.
\]
Since $g$ is zero on the diagonal of $S^2$ this quantity can also be expressed using $P^{(-)}$, and then, since the spectral radius of 
$P^{(-)}$ is less than $\tfrac{1}{2}$ by hypothesis, we obtain the desired convergence to zero as $n$ tends to infinity.
\end{proof}
 
\section{Dynamics and non-endogeny}

Let $L$ be the operator on $L^2(S,{\cal S}, \mu)$ associated with a matrix $\bigl( L(x,x^\prime); x,x^\prime \in S\bigr)$,
\[
Lf(x)= \sum_{x^\prime \in S} L(x,x^\prime) f(x^\prime).
\]
We use the transition probabilities $P^{(2)}$ to determine a new operator ${\cal P}L$  also  acting on $L^2(S,{\cal S},\mu)$, which is  
defined via its associated matrix  via,
\begin{equation}
\label{calp}
(\mu \otimes L ) P^{(2)} = \mu \otimes ({\cal P} L).
\end{equation}
Here  $(\mu \otimes L) (x,x^\prime) = \mu(x)L(x,x^\prime)$ is treated as a row vector on $S^2$ and $P^{(2)}$ acts by matrix multiplication 
on 
the right.
In a similar way we may define a quadratic superoperator ${\cal Q}$ by using the mapping $T^{(2)}$,
\begin{equation}
\label{calq}
T^{(2)}(\mu \otimes L, \mu \otimes L) = \mu \otimes ({\cal Q} L).
\end{equation}
Both superoperators ${\cal P}$ and ${\cal Q}$ arise when considering the isometric embedding of Hilbert spaces:
$ L^2(S,{\cal S}, \mu) \rightarrow L^2(S^2\times E, {\cal S}^2\times {\cal E}, \otimes^2 \mu \otimes \nu)$
 given by $ f \mapsto  f \circ \phi$.  Given  $L$ acting on   on $L^2(S,{\cal S}, \mu)$, it is easily verified that the  
new operator ${\cal P}L$ satisfies
\begin{equation}
\label{forms}
\bigl( f\circ \phi , (L \otimes I \otimes I) g\circ \phi \bigr)_{L^2(S^2\times E, {\cal S}^2\times {\cal E}, \otimes^2 \mu \otimes \nu)}= 
\bigl( f, ({\cal P} L) 
g\bigr)_{  L^2(S,{\cal S}, \mu)},
\end{equation}
for all $f,g \in  L^2(S,{\cal S}, \mu)$.
Similarly  $L$ and ${\cal Q}L$ satisfy
\begin{equation}
\label{forms2}
\bigl( f\circ \phi , (L \otimes L \otimes I) g\circ \phi \bigr)_{L^2(S^2\times E, {\cal S}^2\times {\cal E}, \otimes^2 \mu \otimes \nu)}= 
\bigl( f, ({\cal Q} L) 
g\bigr)_{  L^2(S,{\cal S}, \mu)},
\end{equation}
for all $f,g \in  L^2(S,{\cal S}, \mu)$.
It is a  consequence of the recursive structure  that these two relations extend  as is recorded in the following proposition whose proof we 
omit.
\begin{proposition}
\label{levels}
If $v \in G_{n+1}$ is a daughter of some $u \in G_n$ then, for all $f,g \in {\cal H}_n$,
\[
\bigl( f,L^{(v)} g\bigr)= \bigl( f, ({\cal P} L)^{(u)} g \bigr).
\]
If $v_1 ,v_2 \in G_{n+1}$ are daughters of distinct $u_1,u_2 \in G_n$, then for $f,g \in {\cal H}_n$,
\[
\bigl( f, L^{(v_1)} L^{(v_2)} g\bigr)= \bigl( f, ({\cal P} L)^{(u_1)}({\cal P} L)^{(u_2)} g \bigr).
\]
If $v_1 ,v_2 \in G_{n+1}$ are the two daughters of some  $u \in G_n$, then for $f,g \in {\cal H}_n$,
\[
\bigl( f, L^{(v_1)} L^{(v_2)} g\bigr)= \bigl( f, ({\cal Q} L)^{(u)} g \bigr).
\]
\end{proposition}

Corresponding to  the  principal  eigenvalue $\rho$ of $P^{(-)}$ is a  left eigenvector $\kappa$  satisfying $\kappa P^{(-)} =\rho \kappa$.  
Of course $\kappa$ is only determined up to a scalar multiple and we make some  arbitrary  choice.  Considered as a function on $S \setminus 
S^\nearrow$, $\kappa$ is symmetric since $P^{(-)}$ preserves the space of vectors with this symmetry. It is conceivable that $P^{(-)}$ is 
not irreducible in which case there may be some  further freedom in choosing $\kappa$. This does matter so long as we always choose it, as 
we may, to be symmetric. Next we define a symmetric operator $Q$ on $L^2(S,{\cal S}, \mu)$ from $\kappa$ via
\begin{equation}
\label{q}
\bigl( f, Q g\bigr)_{L^2(S,{\cal S}, \mu)}=-\frac{1}{2}\sum_{(x,x^\prime) \in S^2\setminus S^\nearrow} 
\bigl(f(x^\prime)-f(x)\bigr)\bigl(g(x^\prime)-g(x)\bigr)\kappa(x,x^\prime).
\end{equation}
$Q$ is the generator of an  $S$-valued Markov process which jumps from $x$ to $x^\prime$ at rate $Q(x,x^\prime)=\kappa(x,x^\prime)/ 
\mu(x)$.  We can always assume that $\mu(x)>0$ for all $ x \in S$ by deleting part of $S$ if necessary.

Define an operator on ${\cal H}_n$ via
\begin{equation}
Q_n=( 2\rho)^{-n}  \sum_{u \in G_n} Q^{(u)}.
\end{equation}
 This operator is the generator of a Markov process taking values in the finite configuration space
$\Omega_n$.
The coordinates $\xi_u$, with $u \in G_n$,  evolve independently, each a copy of the process generated by $Q$
but  with their speed altered by the  factor $(2\rho)^{-n}$. At any instant the $\xi_u$ coordinates for $u \in \Gamma_{n-1}$ are 
determined  from the $\xi_u$ 
coordinates with $u\in G_n$
by application of the map $\phi$ with the innovations $(\epsilon_u; u \in \Gamma_{n-1})$ fixed for all time.

The significance of the family of generators $Q_n$ is that they have a certain consistency property that manifests itself at the level of 
the corresponding  forms. We introduce  the forms  ${\cal E}_n$   defined  ${\cal H}_n$ via
\begin{equation}
{\cal E}_n\bigl( f,  g\bigr)= -\bigl( f, Q_n g\bigr) \qquad \text{ for } f,g \in {\cal H}_n.
\end{equation}
\begin{lemma}
\label{qq} 
The operator $Q$ satisfies 
$
{\cal P} Q= \rho Q,
$
and consequently the forms ${\cal E}_n$ are consistent in the sense that for any $m \leq n$,
\begin{equation*}
{\cal E}_m(f, g) = {\cal E}_n(f, g) \text{ for all } f,g \in {\cal H}_m.
\end{equation*}
\end{lemma}
\begin{proof}
In view of the relation between ${\cal P}$ and $P^{(2)}$ given by \eqref{calp}, to prove the first assertion we must verify that
\[
(\mu \otimes Q) P^{(2)}= \rho (\mu \otimes Q).
\]
$(\mu \otimes Q)(x,x^\prime)= \kappa(x,x^\prime)$ for $x \neq x^\prime$ and so the desired equality holds on $S\setminus S^\nearrow$ since 
there it becomes $\kappa P^{(-)}= \rho \kappa$. We deduce that the equality must also hold on the diagonal by observing that ${\cal P}$ 
preserves the class of operators satisfying $L 1=0$.

Suppose that $f,g \in {\cal H}_m$ then using what we have just shown together with Proposition \ref{levels} we obtain
\begin{multline*}
{\cal E}_{m+1}(f,g)=- (f,Q_{m+1}g)= -( 2\rho)^{-(m+1)}  \sum_{v \in G_{m+1}} \bigl(f,Q^{(v)}g\bigr) \\
=-2( 2\rho)^{-(m+1)}  \sum_{u \in G_{m}} \bigl(f,({\cal P}Q)^{(u)}g\bigr)
=- ( 2\rho)^{-m}  \sum_{u \in G_{m}} \bigl(f,Q^{(u)}g\bigr)={\cal E}_m(f,g),
\end{multline*}
which proves the consistency of the forms.

\end{proof}

It follows from the consistency of the ${\cal E}_n$ just established that  we can define a form ${\cal E}$ on the dense subspace  $\bigcup_n 
{\cal H}_n$ of ${\cal H}$ via  ${\cal E}(f,g)={\cal E}_n(f,g)$ 
whenever $f,g \in {\cal 
H}_n$. However it is not necessarily true that ${\cal E}$ is closable.
\begin{lemma}
If $2\rho>1$ then  for $f \in {\cal H}_m$ for some $m$,
\[
\sup_{n \geq m}  \|Q_n f \| <\infty. 
\]
\end{lemma}
\begin{proof}
Consider $f \in {\cal H}_n$.  Expanding $Q_{n+1}$ as a sum and using  Proposition \ref{levels}, plus ${\cal P} Q= \rho Q$,  gives
\begin{multline*}
\big( f, Q_{n+1}^2 f \bigr)
= (2\rho)^{-2n-2} \sum_{u,v \in G_{n+1}} \bigl(f, Q^{(u)}Q^{(v)} f\bigr)\\
= 2(2\rho)^{-2n-2}\sum_{u \in G_n} \left( f, ({\cal Q} Q +{\cal P} Q^2)^{(u)} )f \right) +
 4(2\rho)^{-2n-2}\sum_{\substack{u,v \in G_n \\ u \neq v}} \left( f,  ({\cal P}Q)^{(u)}({\cal P}Q)^{(v)}f \right)\\
=2(2\rho)^{-2n-2}\sum_{u \in G_n} \left( f, ({\cal Q} Q+{\cal P} Q^2)^{(u)}f \right) - (2\rho)^{-2n}\sum_{u \in G_n } 
\left( f,  (Q^2)^{(u)})f \right) + \bigl(f, Q_n^2 f\bigr)\\
=(2\rho)^{-2n}\sum_{u\in G_n}  \left(f, \hat{Q}^{(u)}f \right) + \bigl(f, Q_n^2 f\bigr),
\end{multline*}
where the operator  $\hat{Q}$ acting on $L^2(S,{\cal S}, \mu)$ is given by
\begin{equation*}
 \hat{Q}= 2(2\rho)^{-2}\bigl({\cal Q}Q +{\cal P}Q^2 \bigr)-  Q^2.
\end{equation*}
But now we compare  the operator $\sum_{u \in G_n} \hat{Q}^{(u)}$ with the number operator $A_n$ defined by \eqref{number}. Notice 
that the constant 
${1}\in L^2(S,{\cal S},\mu)$ satisfies $\hat{Q}{ 1}=0$. Thus each subspace ${\cal H}_S$ of ${\cal H}_n$ is an invariant subspace for 
$\sum_{u \in G_n} \hat{Q}^{(u)}$ whose restriction to ${\cal H}_S$ has norm $|S|\|\hat{Q}\|$. 
Thus  we have
\[
\sum_{u\in G_n}\left( f,   \hat{Q}^{(u)} f \right) \leq \|\hat{Q}\| \bigl(f,A_nf\bigr).
\]
If we express $\bigl(f,A_nf\bigr)$ in terms of $P^{(-)}$ as we did in the previous section   then we  find that it  is bounded by some 
constant 
times $(2\rho)^n$, and  thus   we deduce that, for an appropriate constant $C$,
\[
\left(f,Q^2_{n+1}f \right) \leq \bigl(f,Q^2_n,f\bigr)+C(2\rho)^{-n}.
\]
If $2\rho>1$ then the desired conclusion follows.
\end{proof}

\begin{proof}[Proof of Theorem 1: Supercitical case.]

The consistency of the forms ${\cal E}_n$  can be expressed in terms of the corresponding operators as $Q_mf=P_{{\cal H}_m}Q_nf$ for $f \in 
{\cal H}_m$ and  $n \geq m$. Consequently if $\sup_{n \geq m} \|Q_m f\|< \infty $ then, as $n$ tends to infinity,  $Q_n f$ converges in 
${\cal H}$ to some limit we denote by $\mathring{Q}_\infty f$. We see that $-\mathring{Q}_\infty$ is a positive symmetric operator  with
\[
{\cal E}(f,g)= - \left( f, \mathring{Q}_\infty g \right) \qquad \text{ for } f,g \in \bigcup_n {\cal H}_n,
\]
and consequently ${\cal E}$ is closable [ see Theorem X.23 of \cite{rs2}]. 

By construction the subspace ${\cal K}_n$ lies in the kernel of the operator $Q_n$ for each $n$. Hence
\[ {\cal E}(f,f)= 0 \text { for all  } f \in \bigcup_n {\cal K}_n.
\]
Denoting the closure of ${\cal E}$ by $\bar{\cal E}$  we deduce that 
\[
\bar{{\cal E}} (f,f)=0 \text { for all  } f \in {\cal K}.
\]
But  $\bar{\cal E}(f,f)$ cannot be identically zero on ${\cal H}$ since $Q$, and hence ${\cal E}_0$, was not zero, thus ${\cal K} \neq {\cal 
H}$.
\end{proof}

\section{The critical case}

The critical case $2\rho=1$ is endogenous,  provided we  
impose two additional
non-degeneracy conditions:
\begin{equation}
\label{nondegen1}
{\cal H}_0 \cap {\cal K}^\bot \text{ is trivial;}
\end{equation}
\begin{equation}
\label{nondegen2}
P^{(-)} \text{ is irreducible}.
\end{equation}
The first of these conditions  can thought of as analogous  to the condition on a multitype branching process that every initial condition
 leads to  a non-zero probability of eventual extinction. The  example given in the next paragraph  suggests that it is not possible to 
dispense with some condition of this type. The second condition is  probably not essential, but without it the proof given below would be 
considerably more complicated.

The following example shows how it is possible for  the process to be non-endogenous even if $2 \rho=1$.
 Let $S=\{-1,+1\}$ and $E=\{0,1\}$ with
 $\nu(0)=\nu(1)=1/2$ and $\mu(-1)=\mu(+1)=1/2$. 
 Suppose that
\begin{equation}
\phi(x_0,x_1,z)= {\mathbf 1}(z=0) x_0+ {\mathbf 1}(z=1)x_1.
\end{equation}
$S^2\setminus S^\nearrow=\{ (-1,+1),(+1,-1)\}$ and the transition matrix $P^{(-)}$ is $1/2$ times the identity matrix,
 so plainly $2\rho=1$. Furthermore $\xi_\emptyset \in {\cal H}_0$ is orthogonal to every subspace ${\cal K}_n$ and hence to ${\cal K}$, in 
particular this shows that endogeny does not hold.
However this example does not satisfy the strong symmetry condition we have assumed for $\phi$, namely that 
$\phi(x_0,x_1,z)=\phi(x_1,x_0,z)$ for all $x_0,x_1 \in S$ 
and $z \in E$. I do not know whether there are any examples, with $S$ finite, and  this symmetry assumption upheld, but for which 
\eqref{nondegen1} fails.

The sequence of subspaces $ {\cal H}_0 \cap {\cal K}_n^\bot$ is decreasing, and since ${\cal H}_0$ is finite-dimensional,  
\eqref{nondegen1}
can only hold if there exits some $m$ for which
\begin{equation}
 {\cal H}_0 \cap {\cal K}_m^\bot \text{ is trivial.}
\end{equation}
Now consider the quadratic form $\bigl( f,  P_{{\cal K}_m} f\bigr)$. If the preceding condition holds then this quadratic form 
restricted to  $ f \in {\cal H}_0$ is positive definite, and using again the fact that ${\cal H}_0$ is finite-dimensional we deduce that 
there exists an $\epsilon>0$ such that 
\begin{equation}
\label{smallest}
\bigl( f,  P_{{\cal K}_m} f\bigr) \geq \epsilon (f,f) \text{ for all } f \in {\cal H}_0.
\end{equation}
This extends, see the final paragraph of this section, to 
\begin{equation}
\label{dodgy1}
\bigl( f,  P_{{\cal K}_{n+m}} f\bigr) \geq \epsilon^{|S|} (f,f) \text{ for all } f \in {\cal H}_S,
\end{equation}
 for a subset  $S\subseteq G_n$, where $ n \geq 1$ is arbitrary, whilst $m$ is as above.
Decomposing  $f \in {\cal H}_n$ into its components in the subspaces ${\cal H}_S$ as $S$ varies through subsets of $G_n$, and using 
\eqref{split} below, gives,
\begin{equation}
\label{dodgy2}
\bigl( f,  P_{ {\cal K}_{n+m} }f \bigr)
\geq \sum_{ S \subseteq G_n} \epsilon^{|S|} \bigl( f, P_{{\cal H}_S} f\bigr)=
\bigl( f, P_{{\cal K}_n} f \bigr)+\sum_{\substack{S \subseteq G_n \\ S \neq \emptyset}} \epsilon^{|S|} \bigl( f, P_{{\cal H}_S} f\bigr)
\end{equation}
Consequently for a fixed $f$, 
\begin{equation}
\label{tends}
\sum_{\substack{S \subseteq G_n \\ S \neq \emptyset}} \epsilon^{|S|} \bigl( f, P_{{\cal H}_S} f\bigr) \rightarrow 0, 
\end{equation}
as $n\rightarrow \infty$.

The  criticality assumption  that $2 \rho=1$ implies, in the presence of the additional condition \eqref{nondegen2}, that the  sequence of 
matrices $2^nP_n^{(-)}$ is bounded as $n$ varies. Consequently, for a fixed $f$,
\begin{equation}
\bigl(f,A_n f\bigr)=\sum_{ S \subseteq G_n } |S| \bigl( f, P_{{\cal H}_S} f\bigr) 
\end{equation}
is also bounded. The only way that this is consistent with \eqref{tends} is for
\begin{equation}
\sum_{ \substack{S \subseteq G_n \\ S \neq \emptyset}}  \bigl( f, P_{{\cal H}_S} f\bigr) \rightarrow 0,
\end{equation}
which proves endogeny.

A couple of the steps used above need amplification.  Start by considering a generalization of \eqref{iso} and \eqref{subsetdec}. By 
decomposing the tree at level $n$ we obtain a natural isomorphism
\begin{equation}
{\cal H}_{n+m} \cong  \sideset{}{^{G_n}} \bigotimes L^2(\Omega_m, {\cal F}_m,  m_m)  \sideset{}{^{\Gamma_{n-1}}}\bigotimes L^2(E,{\cal 
E}, \nu).
\end{equation}
Splitting $ L^2(\Omega_m, {\cal F}_m,  m_m)$ into the space of constants together its orthogonal complement
$L^2_0(\Omega_m, {\cal F}_m,  m_m)$,
 we obtain the decomposition
\begin{equation}
{\cal H}_{n+m} = \bigoplus_{S \subseteq G_n} {\cal H}_S(m),
\end{equation}
where for each subset $S$ the subspace ${\cal H}_S(m)$ is generated by vectors of the form $\otimes_{u\in G_n} f_u \otimes_{ v \in 
\Gamma_{n-1}} g_v$ with $f_u \in L^2_0(\Omega_m, {\cal F}_m,  m_m)$  for $u \in S$, $f_u$ being a constant for $u \not\in S$, 
and $g_v \in L^2(E,{\cal E}, \nu)$ being unrestricted. 
Notice that for each subset $S$ of $G_n$ the subspace ${\cal H}_S$ is included in ${\cal H}_S(m)$.
There is a corresponding decomposition
\begin{equation}
{\cal K}_{n+m} = \bigoplus_{S \subseteq G_n} {\cal K}_S(m),
\end{equation}
where ${\cal K}_S(m)$  is a subspace of ${\cal H}_S(m)$. The orthogonal projection $P_{{\cal K}_{n+m}}$ acts on  ${\cal H}_{n+m}$ by 
projecting  each  subspace ${\cal H}_S(m)$ onto the  corresponding subspace ${\cal K}_S(m)$. Accordingly if $f \in {\cal H}_{n+m}$ is 
decomposed as $f=\sum f_S$ with $f_S \in {\cal H}_S(m)$ then  
\begin{equation}
\label{split}
 \bigl(f,P_{{\cal K}_{n+m}}f\bigr) = \sum_S \bigl( f_S, P_{{\cal K}_{n+m}} f_S \bigr)
\end{equation}
a fact that was used at \eqref{dodgy2}.

Further examination reveals that ${\cal H}_S(m)$  is naturally isomorphic to a tensor product of ${\cal K}_n$ together with $|S|$ copies 
of ${\cal H}_m^0$, where the latter is the orthogonal complement of the space of constants in ${\cal H}_m$. Similarly ${\cal K}_S(m)$  is 
naturally isomorphic to a tensor product of ${\cal K}_n$ together with $|S|$ copies of ${\cal K}_m^0$.  The restriction of $P_{{\cal 
K}_{n+m}}$ to ${\cal H}_S(m)$ respects this tensor product structure, acting as the identity on the  factor of  ${\cal K}_n$ tensored 
with   copies of the  natural projection from  ${\cal H}_m^0$ to ${\cal K}_m^0$ on the other factors. 
Similarly the restriction of $P_{{\cal H}_n}P_{{\cal 
K}_{n+m}}$ to ${\cal H}_S(m)$ is the tensor product of the identity on ${\cal K}_n$ with $|S|$ copies of  $P_{{\cal H}_0}P_{{\cal K}_m}$ 
acting on ${\cal H}_m^0$. The inequality \eqref{dodgy1} follows from this and \eqref{smallest} since the smallest eigenvalue of a tensor 
product of operators is the product of the smallest eigenvalues.

\section{ Convergence to the dynamics}

Throughout this section we work with the case $2\rho>1$ and we make  the additional assumption
\begin{equation}
 P^{(-)} \text{  is a primitive matrix}.
\end{equation}
According to Perron-Frobenius theory, under this condition,  the limit
of the rescaled $n$-step transition matrices $ \rho^{-n} P_n^{(-)}$ exists and is given by
\begin{equation}
\label{con}
\lim_{n \rightarrow \infty} \rho^{-n}P^{(-)}_n( (x,x^\prime), (y,y^\prime))= \theta(x,x^\prime)\kappa(y,y^\prime),
\end{equation}
where $\theta$ is the left eigenvector of $P^{(-)}$ corresponding to $\rho$, and, as before, $\kappa$ is the right eigenvector. Here we 
normalize $\theta$ and $\kappa$ so that 
\begin{equation}
\sum_{(x,x^\prime) \in S^2\setminus S^\nearrow} \theta( x,x^\prime)\kappa(x,x^\prime)=1.
\end{equation}
We may also normalize so that
\begin{equation}
\sum_{(x,x^\prime)\in S^2\setminus S^\nearrow} \theta(x,x^\prime) \mu(x) \mu(x^\prime)=1,
\end{equation}
which fixes a choice of $\kappa$. We assume throughout this section that $Q$ is defined by \eqref{q} with this choice of $\kappa$.
We deduce that
\begin{equation}
\rho^{-n}( \mu \otimes \mu - \mu^\nearrow) P_n^{(2)} \rightarrow \kappa^\star,
\end{equation}
where $\kappa^\star(x,x^\prime)=\kappa(x,x^\prime)$ when $x\neq x^\prime$ and $\kappa^\star(x,x)= - \sum_{x^\prime \neq x} 
\kappa(x,x^\prime)$. 
Using  the relationship \eqref{calp} between the superoperator ${\cal P}$ and $P^{(2)}$ this may be recast as  
\begin{equation}
\label{con1}
\lim_{n \rightarrow \infty } \rho^{-n} {\cal P}^n P_{{\mathbf 1}^\bot} = -Q.
\end{equation}

Suppose that $f,g  \in {\cal H}_m$ then a straightforward application of Proposition \ref{levels} gives
\begin{equation}
\frac{1}{(2\rho)^n}\bigl ( f, A_n g \bigr) = \frac{1}{2^m\rho^n} \sum_{u \in G_m} \bigl(f,  ({\cal P}^{n-m}P_{{\mathbf 1}^\bot} )^{(u)} g 
\bigr).
\end{equation} 
Thus in view of \eqref{con1} we  deduce that
\begin{equation} \label{con2}
\frac{1}{(2\rho)^n}\bigl ( f, A_n g \bigr) \rightarrow  -\bigl( f, Q_m g \bigr)= {\cal E}(f,g), \qquad \text{ for all } f,g \in {\cal H}_m.
\end{equation} 
Recall that we proved that the form ${\cal E}$ is closable, and let $ Q_\infty$ be the self-adjoint operator associated with
its closure.  If $f \in \bigcup_m {\cal H}_m$ the limit $\mathring{Q}_\infty f=\lim_n Q_n f$ exists and defines a operator 
$\mathring{Q}_\infty$ with domain $\bigcup{\cal H}_m$. The self -adjoint operator $Q_\infty$ is an extension ( the Friedrichs extension) of 
$\mathring{Q}_\infty.$

\begin{proposition} 
\label{con3}
If $f \in {\cal H}_m$ for some $m$, then  as $n$ tends to infinity
\[
(2\rho)^{-n}A_n f \rightarrow  -Q_\infty f,
\]
in the metric topology of ${\cal H}$.
\end{proposition}

\begin{proof}
Convergence in the metric topology is implied by weak convergence together with convergence of the norms. Consequently 
it is sufficient to verify that if $f \in {\cal H}_m$ for some $m$, then as  $n$ tends to infinity,
\[
\bigl\| (2\rho)^{-n}A_n f\bigr\| \rightarrow \bigl\|Q_\infty f\bigr\|,
\]
noting that  weak convergence follows from this and \eqref{con2}.

We begin by computing  $\bigl \| Q_\infty f\bigr\|$.  We have, for $f \in {\cal H}_m$, and $n>m$,
\begin{multline*}
\bigl( f, Q_n^2f \bigr)= (2\rho)^{-2n} \sum_{u \in G_n } \bigl( f, (Q^2)^{(u)} f\bigr)+(2\rho)^{-2n} \sum_{ \substack{u,v \in G_n \\ u \neq 
v}} \bigl( f,  Q^{(u)}Q^{(v)}f \bigr)  \\
= (2\rho)^{-2n} \sum_{u \in G_n } \bigl( f, (Q^2)^{(u)} f\bigr)+2 \sum_{r=m}^{n-1} 2^{2(n-r-1)}(2\rho)^{-2n}\sum_{u \in G_r} \bigl( f, 
({\cal 
Q}{\cal P}^{n-r-1} Q)^{(u)} f \bigr) \\
  + 2^{2(n-m)}(2\rho)^{-2n} \sum_{ \substack{u,v \in G_m \\ u \neq v}} \bigl( f, ( {\cal P}^{n-m} Q)^{(u)}( {\cal P}^{n-m} Q)^{(v)}f \bigr) 
\\
  =(2\rho)^{-2n} \sum_{u \in G_n } \bigl( f, (Q^2)^{(u)} f\bigr)+ 2\sum_{r=m}^{n-1}(2\rho)^{-2(r+1)} \sum_{u \in G_r} \bigl( f, ({\cal Q} 
Q)^{(u)} f \bigr) \\
  + (2\rho)^{-2m} \sum_{ \substack{u,v \in G_m \\ u \neq v}} \bigl( f,  Q^{(u)} Q^{(v)}f \bigr).
\end{multline*}
Now as $n$ tends to infinity $Q_n f$ converges to $Q_\infty f$ in the metric topology of ${\cal H}$, and hence the limit of the lefthandside 
above is $\bigl\| Q_\infty f\bigr\|^2$. Turning to the righthandside the first term tends to zero, and consequently we deduce that
\begin{equation*}
\bigl\|Q_\infty f \bigr\|^2=(2\rho)^{-2m} \sum_{ \substack{u,v \in G_m \\ u \neq v}} \bigl( f,  Q^{(u)} Q^{(v)}f \bigr) +
2 \sum_{r=m}^{\infty}(2\rho)^{-2r-2} \sum_{u \in G_r} \bigl( f, ({\cal Q} Q)^{(u)} f \bigr).
\end{equation*}
A similar calculation is valid for $\|A_n f\|^2$. If we denote $ P_{{\mathbf 1}^\bot}$ by $L$ and $\rho^{-r}{\cal P}^r P_{{\mathbf 1}^\bot}$ 
by $L_r$, then 
\begin{multline*}
(2\rho)^{-2n} \| A_n f\|^2={(2\rho)^{-2n}} \sum_{u \in G_n} \bigl( f, L^{(u)} f \bigr)+ 2 \sum_{r=m}^{n-1} (2 \rho)^{-2r-2} \sum_{u \in G_r} 
\bigl( f, ( {\cal Q}L_{n-r-1})^{(u)} f \bigr ) \\ +(2\rho)^{-2m} \sum_{ \substack{u,v \in G_m \\ u \neq v}} \bigl( f,  L_{n-m}^{(u)} 
L_{n-m}^{(v)}f \bigr). 
\end{multline*}
 Letting $n$ tend to infinity we observe that since $L_n$ tends to $-Q$ we obtain termwise convergence to the expression for $\bigl \| 
Q_\infty  f \bigr \|^2$. To complete the proof we appeal to dominated convergence noting that since ${\cal Q} L_n 1=0,$ we have the estimate
 \[
 \sum_{u \in G_r} \bigl( f, ( {\cal Q}L_{n-r-1})^{(u)} f \bigr ) \leq K \bigl( f, A_r f\bigr),
 \]
 where $K= \sup_n \bigl \| {\cal Q} L_n \bigr \| <\infty$.
\end{proof}
\begin{lemma}
For $f \in {\cal H}_n$,
\[
\|Q_\infty f\| \leq \frac{K}{(2\rho)^{n}} \|A_n f\|,
\]
where $K= \sup\bigl \{ \|Q_\infty f\|:  f \in {\cal H}_0  \text{ and } \|f\|=1\bigr\}$.
\end{lemma}
\begin{proof}
Recall 
the natural isomorphism
\begin{equation*}
{\cal H}_n \cong  \sideset{}{^{G_n}} \bigotimes L^2(S, {\cal S},  \mu)  \sideset{}{^{\Gamma_{n-1}}}\bigotimes L^2(E,{\cal 
E}, \nu).
\end{equation*}
Similarly
\begin{equation*}
{\cal H} \cong  \sideset{}{^{G_n}} \bigotimes L^2(\Omega, {\cal F},  m)  \sideset{}{^{\Gamma_{n-1}}}\bigotimes L^2(E,{\cal 
E}, \nu).
\end{equation*}
The   operator $Q_\infty:{\cal H}_n \rightarrow {\cal H}$  can be written as a sum
\[Q_\infty = \frac{1}{(2\rho)^n}\sum_{u \in G_n} Q_\infty^{(u)},
\]
where $Q_\infty^{(u)}$  is unitary equivalent the tensor product of  
$Q_\infty: L^2(S, {\cal S},  \mu) \rightarrow  L^2(\Omega, {\cal F},  m)$ on the factor corresponding 
to the node $u \in G_n$ and the identity on all other factors. It is easy to see from this structure, together with the fact that
$Q_\infty 1=0$ that 
if $f \in {\cal H}_n$  decomposes as $f =\sum_{S \subseteq  G_n} f_S$ with $f_S\in {\cal H}_S$ then 
\begin{multline*}
\| Q_\infty f\|^2= 
\sum_S  \|Q_\infty f_S\|^2 = \frac{1}{(2\rho)^{2n}} \sum_S \sum_{u,v \in S}  \bigl(Q_\infty^{(u)} f_S, Q_\infty^{(v)} f_S\bigr) \\
\leq \frac{K^2}{(2\rho)^{2n}} \sum_S  |S|^2 \|f_S\|^2 =  \frac{K^2}{(2\rho)^{2n}} \|A_n f\|^2.
\end{multline*}
\end{proof}

\begin{lemma}
$\bigcup_m {\cal H}_m$ is a core for $Q_\infty$.
\end{lemma}
\begin{proof}
Recall that $\mathring{Q}_\infty$ is the restriction of $Q_\infty$ to $\bigcup {\cal H}_m$. 
To verify the claimed result   it suffices to show that, for some $\alpha>0$,  the range of $(\alpha-\mathring{Q}_\infty)$  is dense in 
${\cal H}$.
 For this shows that $\alpha$ belongs to the resolvent set of the closure of $\mathring{Q}_\infty$, and then we apply the 
criterion 
of Theorem X.1 of \cite{rs2}. 

Fix some  $\alpha>0$. Let $R^\alpha_n: {\cal H}_n \rightarrow {\cal H}_n$ be  the $\alpha$ resolvent of $Q_n$.
Given $f\in \bigcup_m {\cal H}_m$ let $v_n= R^\alpha_n f$ for $n$ sufficiently large.   Let $g\in \bigcup {\cal H}_m$, then for  all 
sufficiently large $n$,
\[
\bigl( g, (\alpha-Q_\infty)v_n\bigr)=   \bigl( g, (\alpha-Q_n)v_n\bigr) =(g,f).
\]
Now suppose that we know that $(\alpha-Q_\infty)v_n$ is uniformly bounded in norm, then we deduce that, for any $g \in {\cal H}$,
\[
\bigl( g, (\alpha-Q_\infty)v_n\bigr) \rightarrow (g,f) \text{ as } n \rightarrow \infty.
\]
Thus the range of $(\alpha-\mathring{Q}_\infty)$  is weakly dense, and consequently norm dense in ${\cal H}$.

To verify the supposition that $(\alpha-Q_\infty)v_n$ is uniformly bounded in norm we note that since $v_n$ are uniformly
  bounded in suffices to verify that $Q_\infty v_n$ are also.  Using the previous lemma  and the fact that $R_n^\alpha$ and $A_n$ commute 
we have
\begin{multline*}
\bigl\| Q_\infty v_n \bigr \| \leq  \frac{K}{(2\rho)^{n}} \bigl\| A_n v_n\bigr\|
=\frac{K }{(2\rho)^{n}}\bigl\|A_n R_n^\alpha f \bigr\|= \frac{K}{(2\rho)^{n}}\bigl\| R_n^\alpha A_n f \bigr\| 
\leq \frac{K}{\alpha(2\rho)^{n}} \bigl\|  A_n f \bigr\|,
\end{multline*}
and the righthandside is bounded as $n$ tends to infinity.
 \end{proof}

\begin{proof}[Proof of Theorem \ref{thmcon}]
 We know that $\cup {\cal H}_n$ is a common core for $Q_\infty$ and for $ A_n P_{{\cal H}_n}$. Thus the convergence established at 
Proposition 
\ref{con3} implies 
 strong resolvent convergence of $(2\rho)^{-n} A_n P_{{\cal H}_n}$ to $-Q_\infty$  (see  Theorem VIII.25 of \cite{rs}  ).
\end{proof}

Convergence of the generators in the strong resolvent sense implies that the semigroups converge in the strong operator topology.  From this 
fact we 
 obtain the following corollary, which can can also be expressed in terms of $T^{(2)}$. 

\begin{corollary} For each fixed $t>0$, as $n \rightarrow \infty$, 
\[
{\cal Q}^n \bigl( P_{\mathbf  1} + e^{-(2\rho)^{-n} t} P_{{\mathbf 1}^\bot } \bigr)\rightarrow M_t,  
\]
where $M_t$ acting on $L^2(S)$ is defined by identifying  ${\cal H}_0$ with $L^2(S)$ and then setting $ \bigl(f,M_t g\bigr)_{L^2(S)}= \bigl( f , 
e^{tQ_\infty} g \bigr )$.
\end{corollary}
Another interpretation of this corollary is available in terms of the spectral measures  $\mu_f^{(n)}$ defined at \eqref{specmeasure}. The 
rescaled measures
\begin{equation}
\tilde{\mu}_f^{(n)} \bigl ([0,x]\bigr) = \mu_f^{(n)}\bigl(([0, (2\rho)^n x]\bigr) \qquad \text{ for }x \geq 0,
\end{equation}
converge weakly towards a measure $\mu_f$ whose  Laplace transform is given by
\begin{equation}
\int_0^\infty e^{-tx} \mu_f(dx) = \bigl( f,  e^{tQ_\infty} f \bigr).
\end{equation}

\end{document}